\documentclass{amsart}

\usepackage{amsmath,amssymb}
\usepackage{mathtools}
\usepackage{lmodern}

\makeindex

\newtheorem{theorem}{Theorem}

\newtheorem{remark}[theorem]{Remark}

\newcommand{\C}{\mathbb{C}}
\newcommand{\B}{\mathbb{B}}
\newcommand{\id}{\mathrm{id}}

\begin{document}

\title[Quantum invariants]{Quantum invariants via Hopf algebras and solutions to the Yang-Baxter equation}
\author{Leandro Vendramin}
\address{IMAS--CONICET and Depto. de Matem\'atica FCEN,
Universidad de Buenos Aires, Pab. I -- Ciudad Universitaria (1428),
Buenos Aires, Argentina}
\email{lvendramin@dm.uba.ar}

\maketitle

\begin{flushright}
	A Carloncha
\end{flushright}

\medskip
The fundamental problem of knot theory is to know whether two knots are
equivalent or not. As a tool to prove that two knots are different,
mathematicians have developed various invariants.  Knots invariants are just
functions that can be computed from the knot and depend only on the topology of
the knot.  Here we describe quantum invariants, a powerful family of invariants
related to the celebrated Yang--Baxter equation. 

The basic idea of quantum invariants belongs to Witten~\cite{MR990772}. It was
then developed by Reshetikhin and Turaev in~\cite{MR1036112} and
~\cite{MR939474}. In modern language, Reshetikhin and Turaev's construction is
essentially the observation that finite-dimensional representations of a
quantum group $\mathcal{U}_q(\mathfrak{g})$ form a ribbon category. 

\subsection*{Braidings and Yetter-Drinfeld modules}

Recall that a \emph{braided vector space} is a pair $(V,c)$, where $V$ is a
vector space and $c\colon V\otimes V\to V\otimes V$ is a solution of the
\emph{braid equation}, that is a linear isomorphism such that
\[
(c\otimes\id)(\id\otimes c)(c\otimes\id)=
(\id\otimes c)(c\otimes\id)(\id\otimes c).
\]

There is a procedure that produces braided vector spaces from Hopf algebras. We
will present this method in the language of Yetter-Drinfeld modules. Let $H$ be
a Hopf algebra with invertible antipode $S$.  A \emph{Yetter-Drinfeld module}
over $H$ is triple $(V,\cdot,\delta)$, where $(V,\cdot)$ is a left $H$-module,
$(V,\delta)$ is a left $H$-comodule and such that the compatibility
\[
\delta(h\cdot v)=h_1v_{-1}S(h_3)\otimes h_2v_0
\]
holds for all $h\in H$ and $v\in H$. 

Yetter-Drinfeld modules over $H$ form a category (the objects are the
Yetter-Drinfeld modules over $H$ and the morphisms are linear maps that are
module and comodule homomorphisms). If $V$ and $W$ are Yetter-Drinfeld modules
over $H$, then $V\otimes W$ is a Yetter-Drinfeld module over $H$
with
\[
h\cdot v\otimes w=h_1\cdot v\otimes h_2\cdot w,
\quad
\delta(v\otimes w)=v_{-1}w_{-1}\otimes (v_0\otimes w_0).
\]

Yetter-Drinfeld modules produce braided vector spaces. More generally, if $U$,
$V$ and $W$ are Yetter-Drinfeld modules over $H$, the map 
\[
c_{V,W}\colon V\otimes W\to W\otimes V,
\quad
c_{V,W}(v\otimes w)=v_{-1}\cdot w\otimes v_0,
\]
is an invertible morphism of Yetter-Drinfeld modules over $H$ such that 
\[
(c_{W,X}\otimes\id_V)(\id_W\otimes c_{V,X})(c_{V,W}\otimes\id_{X})
=(\id_X\otimes c_{V,W})(c_{V,X}\otimes\id_W)(\id_V\otimes c_{W,X})
\]
In particular, when $X=V=W$ one obtains braided vector spaces. 

\begin{remark}
Yetter-Drinfeld modules and their braided vector spaces appear in the
classification of Hopf algebras, see for example~\cite{MR1994219}.
\end{remark}

Finite-dimensional Yetter-Drinfeld modules have duals. Let $V$ be a
finite-dimensional Yetter-Drinfeld module over $H$ with basis
$\{v_1,\dots,v_n\}$. Let $V^*$ be the vector space with basis
$\{f_1,\dots,f_n\}$, where $f_i(v_j)=\delta_{ij}$. Then $V^*$ is a
Yetter-Drinfeld module over $H$ with
\[
(h\cdot f)(v)=f(S(h)\cdot v),\quad
\delta(f)=\sum_{k=1}^n S^{-1}\left( (f_i)_{-1} \right) \otimes f((v_i)_0)f_i
\]
for all $h\in H$, $f\in V^*$ and $v\in V$.  
Recall that there are two canonical maps $b_V\colon\C\to V\otimes
V^*$, $e_V\colon V^*\otimes V\to\C$, such that 
\[
 	(\id_V\otimes e_V)(b_V\otimes\id_V)=\id_V,\quad
     (e_V\otimes\id_{V^*})(\id_{V^*}\otimes b_V)=\id_{V^*}.
\]
The map $b\colon\C\to V\otimes V^*$ is given by
\[
b(1)=\sum_{j=1}^n v_i\otimes f_i
\]
and $e\colon V^*\otimes V\to\C$, $f\otimes v\mapsto f(v)$, is the usual
evaluation map. We leave to the reader the exercise of proving that
\[
c_{V^{*},V}=(e\otimes\id_V\otimes\id_{V^{*}})(\id_{V^*}\otimes c_{V,V}\otimes\id_{V^*})(\id_{V^*}\otimes\id_V\otimes b).
\]
Similarly one computes $c_{V,V^*}$ and $c_{V^*,V^*}$.


The category of finite-dimensional Yetter-Drinfeld modules over $H$ form a
\emph{ribbon} category, this means that there is a family of natural
isomorphisms $\theta_V\colon V\to V$ such that 
\[
\theta_{\C}=\id_{\C},
\quad
\theta_{V^*}=(\theta_{V})^*,\quad
\theta_{V\otimes W}=c_{W,V}c_{V,W}(\theta_V\otimes\theta_W).
\]
See~\cite{MR1231205} for a classification of all such $\theta$.

\subsection*{Braids and knots}

The \emph{braid group} $\B_n$ is the group with generators
$\sigma_1,\dots,\sigma_{n-1}$ and relations
\begin{align*}
& \sigma_i\sigma_j=\sigma_j\sigma_i && |i-j|\geq2,\\
& \sigma_i\sigma_{i+1}\sigma_i=\sigma_{i+1}\sigma_i\sigma_{i+1} && 1\leq i\leq n-2.
\end{align*}

\begin{figure}[h]
    \includegraphics[scale=0.2]{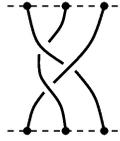}
	\centering
     \caption{The braid $\sigma_1^{-1}\sigma_2\sigma_1^{-1}\in\B_3$.}
     \label{fig:braid}
\end{figure}

One can draw pictures to represent elements of the braid group, see for example
Figure~\ref{fig:braid}.  The group operation is just composition of braids.

Braided vector spaces produce representations of the braid group: if $(V,c)$ is
a braided vector space, the map $\rho_n\colon\B_n\to\mathbf{GL}(V^{\otimes})$, 
$\sigma_j\mapsto c_j$, where
\[
c_j(v_1\otimes\cdots\otimes v_n)
=v_1\otimes\cdots\otimes v_{j-1}\otimes c(v_j\otimes v_{j+1})\otimes v_{j+2}\otimes\cdots\otimes v_n,
\]
extends to a group homomorphism.  This $\rho_n$ yields a
graphical calculus where certain morphisms are represented by pictures.

\begin{figure}[h]
    \includegraphics[scale=0.2]{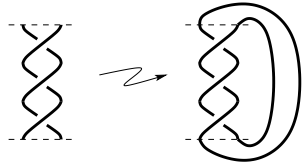}
	\centering
    \caption{The closure of a braid.}
    \label{fig:closure}
\end{figure}

One can use braids to represent knots. Alexander's theorem states that each knot
is the closure of a braid, see Figure~\ref{fig:closure}.  This 
gives a way to move between braids and knots. Two equivalent braids will become
equivalent knots. However, the closure of two non-equivalent braids could
produce equivalent knots. This means that one needs a method to determine if
two closed braids are equivalent. The answer to this fundamental method is now
known as Markov's theorem. 

\begin{figure}[h]
    \includegraphics[scale=0.4]{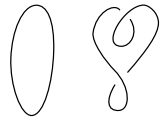}
	\centering
	\caption{Trivial knots.}
	\label{fig:trivial}
\end{figure}

To obtain invariants of knots from quantum groups, one needs to accept that the
invariants on the diagrams of Figure~\ref{fig:trivial} could receive different
values. Without this it would be very hard to construct good invariants. A
\emph{framed knot} is a knot equipped with a smooth family of non-zero vectors
orthogonal to the knot. Once a plane projection is given, one way to choose a
framing is to use a vector field everywhere parallel to the projection plane.
This is the \emph{blackboard framing}. Therefore a framed knot can be viewed as
a \emph{ribbon knot}. Framed knots with blackboard framing are not invariant
under the first Reidemeister move.

We can use morphisms of the category of finite-dimensional Yetter-Drinfeld modules over $H$ 
to color the arcs of our knot. We use $c$,
$c^{-1}$, $e$, $b$,
\begin{align*}
&b^{-}\colon\C\to V^{*}\otimes V,&& b^-=(\id_V\otimes\theta^{-1}_V)c^{-1}_{V,V^*}b,\\
\shortintertext{and}
&e^-\colon V\otimes V^{*}\to\C,  &&
e^-=ec_{V,V^*}(\theta_V\otimes\id_{V^*})
\end{align*}
as follows:
\begin{figure}[h]
\begin{align*}
&\raisebox{-22pt}{\makebox(15,10){$c=$}}\raisebox{-40pt}{\includegraphics[width=65pt]{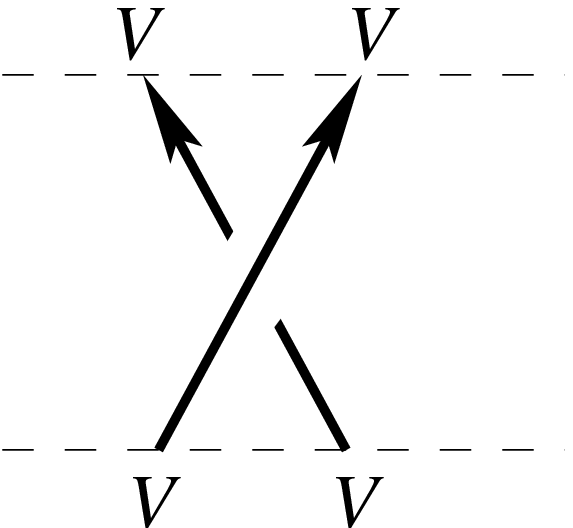}} && 
\raisebox{-22pt}{\makebox(15,10){$c^{-1}=$}}\raisebox{-40pt}{\includegraphics[width=65pt]{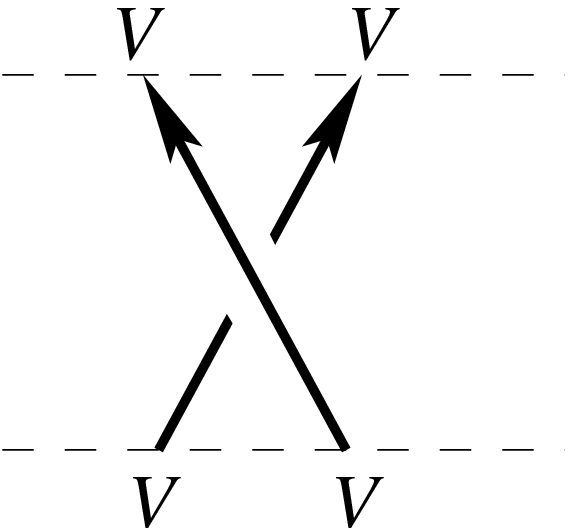}} \\
&\raisebox{-22pt}{\makebox(15,10){$e^-=$}}\raisebox{-40pt}{\includegraphics[width=65pt]{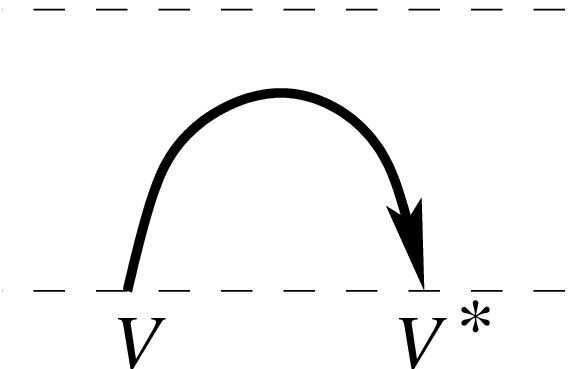}} && 
 \raisebox{-22pt}{\makebox(15,10){$e=$}}\raisebox{-40pt}{\includegraphics[width=65pt]{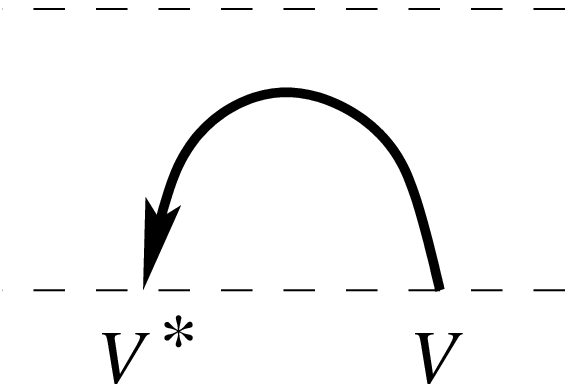}} \\
&\raisebox{-22pt}{\makebox(15,10){$b^-=$}}\raisebox{-40pt}{\includegraphics[width=65pt]{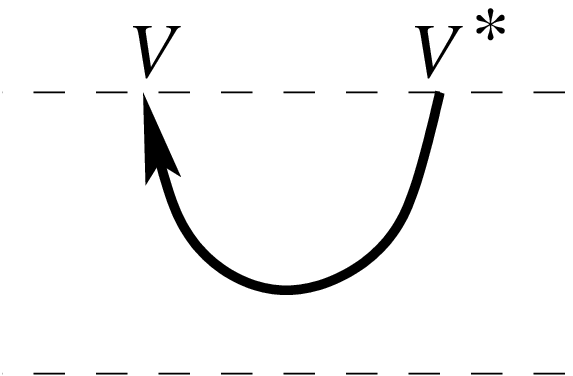}} && 
 \raisebox{-22pt}{\makebox(15,10){$b=$}}\raisebox{-40pt}{\includegraphics[width=65pt]{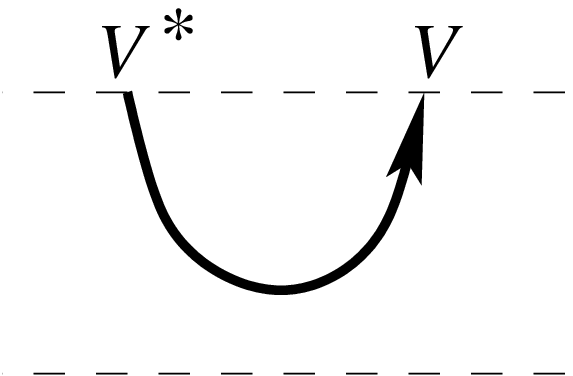}} 
\end{align*}
\end{figure}

It is an exercise to check that this is invariant under the second and third
Reidemeister moves.  In some cases this construction yields a knot invariant.
Let us assume for example that the twist $\theta$ acts on $V$ and $V^*$ by scalar
multiplication. To normalize the diagram one multiplies by
$\theta^{-\omega(K)}$, where $\omega(K)$ is the writhe of the knot $K$ (that is
the number of positive crossings of $K$ minus the number of negative crossings
of $K$). 

\subsection*{An example: The Taft algebra}

Let $m$ be an odd positive integer and let $q\in\C^\times$ be a primitive
$m$-root of $1$. Let $H$ be the algebra generated by $x$ and $g$ with relations
\[
gx=qxg,\quad
x^m=0,\quad
g^m=1.
\]
Then $\dim H=m^2$ and $\{x^ig^j:1\leq i,j<m\}$ is a basis of $H$. 
Furthermore $H$ is a Hopf algebra with
\begin{align*}
&\Delta(g)=g\otimes g,
&&\Delta(x)=g\otimes x+x\otimes 1,
&&\epsilon(g)=1,\\
&\epsilon(x)=0,
&&S(g)=g^{-1},&&S(x)=-g^{-1}x.
\end{align*}

Let $n\geq0$ and $V_n$ be a vector space with basis
$\{v_{-n},v_{-n+2},\dots,v_{n-2},v_n\}$. Then $V_n$ is a simple left $H$-module
with 
\[
g\cdot v_k=q^{-k/2},\quad
x\cdot v_k=v_{k-2},
\]
and a left $H$-comodule with
\begin{gather*}
\delta(v_k)=\sum_{i\geq0}\frac{1}{(i)!_q}\prod_{j=0}^{i-1}\alpha_{k+2j}x^ig^{-\frac{k+2i}{2}}\otimes v_{k+2i},
\shortintertext{where}
\alpha_k=q^{-\frac{k+n+2}{2}}(q-1)\left(\frac{k+n+2}{2}\right)_q\left(\frac{n-k}{2}\right)_q
\end{gather*}
and $(k)_q=1+q+\cdots +q^{k-1}$.

One proves that the $V_n$ are simple Yetter-Drinfeld modules over $H$. The
simplicity of $V_n$ implies that $\theta$ acts on $V_n$ by scalar multiplication by some
number $\theta_n$. Since 
\[
c_{V_n,V_m}(v_n\otimes v_m)=q^{nm/4}v_m\otimes v_n,
\]
one needs $\theta_{n+m}=q^{nm/2}\theta_n\theta_m$ and $\theta_0=1$. By
using that $V_1\otimes V_1\simeq V_2\oplus V_0$, one proves that
$\theta_n=q^{\frac{n^2+2n}{4}}$. 

Let us now study the case $n=1$. The matrix of the braiding is 
\[
\begin{pmatrix}
	q^{1/4} & 0 & 0 & 0\\
	0 & q^{1/4}-q^{-3/4} & q^{-1/4} & 0\\
	0 & q^{-1/4} & 0 & 0\\
	0 & 0 & 0 & q^{1/4}
\end{pmatrix}.
\]
We leave to the reader the exercise to compute the matrices of $c_{V^*,V^*}$,
$c_{V^*,V}$ and  $c_{V^*,V}$. Then
\begin{align*}
&b^{-}(1)=q^{1/2} v_{-1}\otimes v_{-1}+q^{-1/2}v_1\otimes v_1,\\
&e^{-}(v_{-1}\otimes v_{-1},v_{-1}\otimes v_1,v_1\otimes v_{-1},v_1\otimes v_1)=(q^{1/2},0,0,q^{1/2}). 
\end{align*}

\begin{figure}[h]
\includegraphics[scale=0.2]{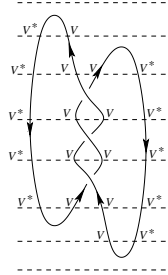}
\centering
\caption{The left trefoil knot.}
\label{fig:trefoil}
\end{figure}

Let us apply this to prove that the left and right trefoil knots are different. These knots
are the closure of $\sigma_1^{-3}\in\B_2$ and $\sigma_1^{3}\in\B_2$.  The
invariant for the left trefoil knot of Figure~\ref{fig:trefoil} is 
\begin{align*}
	\theta_2^{-3}(e\otimes e^-)(\id_{V^*}\otimes c_{V,V}^{-3}\otimes\id_{V^*})(b^-\otimes b)
&=-q^{9/2}+q^{5/2}+q^{3/2}+q^{1/2}.
\end{align*}
The invariant for the right 
trefoil knot is
\[
\theta_2^{3}(e\otimes e^-)(\id_{V^*}\otimes c_{V,V}^{3}\otimes\id_{V^*})(b^-\otimes b)=
q^{-1/2}+q^{-3/2}+q^{-5/2}+q^{1/2}.
\]
Therefore that these knots are different.

\begin{remark}
The modules $V_n$ can be used to recover the Jones polynomial. 
\end{remark}

\begin{remark}
In~\cite{MR1918807} Gra\~na showed that the quandle cocycle invariants in the
sense of~\cite{MR1990571} are quantum invariants related to Yetter-Drinfeld
modules over group algebras. 
\end{remark}

\subsection*{Acknowledgements}

This paper is based on a minicourse given by Mat\'ias Gra\~na in C\'orboba in
2003.  Thanks to Sergei Chmutov, the figures were taken from~\cite{MR2962302}.
The author is partially supported by PICT-201-0147 and MATH-AmSud 17MATH-01.

\bibliographystyle{abbrv}
\bibliography{refs}

\end{document}